\documentclass{amsart}
\usepackage{fullpage}
\usepackage{graphicx}
\usepackage{enumerate}
\input pictex.sty

\newtheorem{thm}{Theorem}[section]

\theoremstyle{definition}

\theoremstyle{remark}

\numberwithin{equation}{section}

\begin{document}

\title{A Riemann-Hilbert problem for skew-orthogonal polynomials}%
\author{Virgil U. Pierce}%
\address{Dept. of Math, The Ohio State University}%
\email{vpierce@math.ohio-state.edu}%

\thanks{This work was 
supported in part by NSF grant DMS-0135308.}%
\subjclass{}%
\keywords{}%

%\date{}%
%\dedicatory{}%
%\commby{}%
% ----------------------------------------------------------------
\begin{abstract}
We find a local $(d+1) \times (d+1)$ Riemann-Hilbert problem characterizing the
skew-orthogonal polynomials associated to the partition function of
the Gaussian Orthogonal Ensemble of random matrices with a potential
function of degree $d$.  Our Riemann-Hilbert problem is
similar to a local  $d \times d$  Riemann-Hilbert problem found by Kuijlaars
and McLaughlin 
characterizing the bi-orthogonal polynomials.  This gives more
motivation for finding methods to compute asymptotics of high order
Riemann-Hilbert problems, and brings us closer to finding asymptotics
of the skew-orthogonal polynomials.
\smallskip

\begin{flushleft}
{\em Keywords: Skew-orthogonal polynomials; Riemann-Hilbert problems}
\end{flushleft}
\end{abstract}
\maketitle
%\tableofcontents
% ----------------------------------------------------------------

\section{Introduction}

The partition function of random matrices for the Gaussian Unitary
Ensemble (GUE) is associated with a family of 
orthogonal polynomials $P_{j}(x)$ 
given by the conditions
\begin{equation} \label{orthopoly}
\int_{\mathbb{R}}
 P_j(x) P_k(x) e^{- N V(x)} dx = h_j \delta_{jk},
\end{equation}
where $V(x)$ is an even degree polynomial with positive leading coefficient,
and $ P_j(x) = x^j + \mathcal{O}(x^{j-1})$ for $x\to \infty$ 
\cite{BIZ, EM, Mehta2} .  
Define the orthonormal polynomials $\hat{P}_j(x) = P_j(x) / h_j^{1/2}$.  
The asymptotic expansion of the classical families of orthogonal
polynomials were computed by \cite{Szego}.
The orthonormal polynomials satisfy a three term recursion relation
that, under changes in $V$, satisfies the Toda lattice hierarchy.
They may be represented as Hankel determinants \cite{Robert}.

To compute the orthogonal polynomials numerically it is helpful to use
the following result:
Define the moments
\begin{equation*}
m_{i} = \int_{\mathbb{R}} x^i e^{- NV(x)} dx .
\end{equation*}
Let $M$ be the symmetric matrix $M_{ij} = m_{i+j}$.  
Then factor $M = Q D R$ using the
$QR$-algorithm.  The columns of $R^{-1}$ give the coefficients of the
orthogonal polynomials.  

One may show that the zeros of the orthogonal
polynomials are real, interlace, and their distribution converges 
to a measure $\mu(V)$
which is given as the solution to a minimization problem
\cite{Deift}.  

The asymptotic expansions of $P_N$ and $P_{N-1}$ 
are computed by realizing them as the unique
solutions of a $2\times 2$ Riemann-Hilbert problem (see \cite{Deift, DKMVZ, DKMVZ2, KMvAV} ).
 These
 expansions were assembled and used to prove the existence of an
 asymptotic expansion of $\frac{1}{N^2} \log\left( Z_N(V) \right)$ for
 a cone of $V$, in \cite{EM}.

Similar studies have been started regarding multi-matrix models and
multi-orthogonal polynomials \cite{Bleher_Kuijlaars}.  The simplest
case is that of  bi-orthogonal 
polynomials   given
by two families of polynomials
$\left\{ P_k(x) \right\}$ and $\left\{ Q_k(x) \right\}$ satisfying 
\begin{equation*}
\iint_{\mathbb{R}^2}  
P_j(x) Q_k(y) \exp\left[ - N ( V(x) + W(y)) -
  2N \tau xy \right] dx dy = h_j \delta_{jk},
\end{equation*}
where $\tau$ is a non-zero coupling constant,
$V(x)$ and $W(y)$ are polynomials of even degree with positive leading
coefficient, 
 $P_j(x) = x^j + \mathcal{O}(x^{j-1})$ for $x\to \infty$, 
and $Q_k(y) = y^k + \mathcal{O}(y^{k-1})$
 for $y \to \infty$
\cite{Ercolani_McLaughlin}.  
The bi-orthogonal polynomials satisfy a recursion relation that, under
 changes in $V$, satisfies the Full-Kostant Toda lattice hierarchy.
One may compute them as determinants \cite{Ercolani_McLaughlin}.
They may be characterized in terms of Riemann-Hilbert problems (see
 \cite{BEH, Ercolani_McLaughlin, Kapaev, KM, AGK}).

To compute the bi-orthogonal polynomials numerically it is helpful to
use the following result:
Let $M$ be the matrix of moments 
\begin{equation*}
M_{ij} = \iint_{\mathbb{R}^2}  x^i y^j \exp\left[ 
- N (V(x) + W(y)) - 2N\tau xy \right] dx dy.
\end{equation*}
 Then factor $M=L D U$ using Gaussian
elimination.  The columns of $U^{-1}$ give the coefficients of
$Q_k(y)$ while the rows of $L^{-1}$ give the coefficients of
$P_j(x)$.  

One may show that the zeros of the bi-orthogonal polynomials are real
\cite{Ercolani_McLaughlin}.  One finds numerically that they interlace
and the suspicion is that their distributions
 converge to a pair of measures.

The bi-orthogonal polynomials are solutions of a non-local
 $2\times 2$ Riemann-Hilbert problem (see \cite{BEH, Ercolani_McLaughlin}).
 Recently a local $d\times d$ Riemann-Hilbert problem that determines the
 bi-orthogonal polynomials was found in the
 case when $\mbox{deg}(W) = d$ by \cite{KM}.  
The advantage of this
 local problem is that there is reason to believe it is more amenable
 to asymptotic analysis, see \cite{McL06}.   
Another high order Riemann-Hilbert problem characterizing the bi-orthogonal
 polynomials has been found by \cite{Bertola}.

Define the skew-inner products 
\begin{equation} 
\label{GOE}
\langle f, g \rangle_1 = \iint_{\mathbb{R}^2} f(x) g(y)
\epsilon(x-y) e^{ - N V(x) - N V(y)} dx dy ,
\end{equation}
where 
\begin{equation*}
\epsilon(x) = \left\{ \begin{matrix}  -1 & : x < 0 \\
1 &: x > 0 
\end{matrix} \right. ;
\end{equation*}
\begin{equation}
\label{GSE}
\langle f, g \rangle_4 = \int_{\mathbb{R}} 
\left[ f(x) g'(x) - f'(x) g(x) \right] e^{ - N V(x) } dx.
\end{equation}

The partition functions of random matrices for the Gaussian Orthogonal
Ensemble (GOE) ($\beta=1$)  or Gaussian Symplectic Ensemble (GSE)
($\beta=4$)  are associated
with skew-orthogonal polynomials $\left\{ p_j(x) \right\}$
given by the conditions
\begin{align} \label{skewortho}
\langle p_{2k}(x), p_{2j}(y) \rangle_\beta
&= 0 \\ 
\langle p_{2k}(x), p_{2j+1}(y) \rangle_\beta
&= h_j \delta_{kj} \\ 
\langle p_{2k+1}(x), p_{2j}(y) \rangle_\beta
& = - h_j \delta_{kj} \\ 
\langle p_{2k+1}(x), p_{2j+1}(y) \rangle_\beta
&= 0 , \label{skewortho_last}
\end{align}
where
$p_{2j}(x) = x^{2j} + \mathcal{O}(x^{2j-1}) $, 
and $p_{2j+1}(x) = x^{2j+1} + \mathcal{O}(x^{2j-1}) $ for $x\to \infty$ 
(see \cite{AFNvM, Deift_Gioev, Deift_Gioev2, Ghosh, Ghosh_Pandey, Mehta, Mehta2,  Nago_Forrester, Stojanovic, Stojanovic_errata,
  Tracy_Widom, Tracy_Widom2, Widom}) .
There is an ambiguity in conditions (\ref{skewortho})-(\ref{skewortho_last}) in the sense that
$\tilde{p}_{2j+1} = p_{2j+1} + a p_{2j}$ may properly be called a
skew-orthogonal polynomial of degree $2j+1$ for any $a$.  We define the
$2j+1$ skew-orthogonal polynomial to be the one with a zero coefficient of
$x^{2j}$.  Define the skew-orthonormal polynomials to be 
$\hat{p}_{2j}(x) = p_{2j}(x) / \langle p_{2j}, p_{2j+1} \rangle_\beta^{1/2},$
and
$ \hat{p}_{2j+1}(x) = p_{2j+1}(x) / \langle p_{2j}, p_{2j+1}
\rangle_\beta^{1/2}$.

The skew-orthonormal polynomials satisfy a recursion relation that,
under changes in $V$, satisfies the Pfaff lattice hierarchy.  For example if
$V = V_0(x) + t_j x^j$ then it can be shown that the recursion
relation is of the form
\begin{equation*}
x \vec{p}(x) = L \vec{p}(x),
\end{equation*}
where $\vec{p}(x) = ( \hat{p}_0(x), \hat{p}_1(x), \hat{p}_2(x), \dots )^T$ 
and 
\begin{equation*}
L = \begin{pmatrix} 
0   & 1   & 0    & 0   &  \\
b_1 & d_1 & a_1  & 0   &  \\
c_1 & e_1 & -d_1 & 1   &  \\
f_1 & c_2 & b_2  & d_2 &  \\
 &  &  &  & \ddots
\end{pmatrix}.
\end{equation*}
Let 
\begin{equation} \label{J}  
J = \begin{pmatrix}
0  & 1  & 0  & 0  &  \\
-1 &  0 & 0  & 0  &  \\
0  &  0 & 0  & 1  &  \\
0  &  0 & -1 & 0  &  \\
 & & & & \ddots
\end{pmatrix}.
\end{equation}
Then $L$ satisfies the differential equation
\begin{equation*}
\frac{dL}{dt_{j}} = - \left[ \pi_{k}( L^j), L \right],
\end{equation*}
where 
\begin{equation*}
\pi_k(M) = M_- - J M_+^T J + \frac{1}{2} ( M_0 - J M_0^T J ) ,
\end{equation*}
and $M_\pm$ is projection onto the upper ( resp. lower)
 $2\times 2$ block triangular parts
of $M$ and $M_0$ is projection onto the diagonal $2\times 2$ blocks of $M$
\cite{AFNvM, AHvM, Adler_Moerbeke} .

Define the Pfaffian of a skew-symmetric matrix $M$ by 
\begin{equation*}
\mbox{det}(M) = \left[ \mbox{pf}(M) \right]^2.
\end{equation*}
Define the skew-symmetric matrix of moments 
\begin{equation*}
M_{ij} = \langle x^i, y^j \rangle_\beta.
\end{equation*}
Then one finds Pfaffian formulas for the skew-orthogonal polynomials
(see \cite{AFNvM, AHvM,  Adler_Moerbeke}):
\begin{equation*}
p_{2j}(x) = 
\mbox{pf}\left[ 
\begin{pmatrix}
0         & M_{01}    & M_{02}    & \dots & M_{0 2j} & 1 \\
-M_{01}   & 0         & M_{12}    & \dots & M_{1 2j} & x \\
\vdots    & \vdots    & \vdots    &       & \vdots   & \vdots \\
-M_{0 2j} & -M_{1 2j} & -M_{2 2j} & \dots & 0        & x^{2j} \\
-1        & -x        & -x^2      & \dots & -x^{2j}  & 0 
\end{pmatrix} \right]
\end{equation*}
and
\begin{equation*}
p_{2j+1}(x) = 
\mbox{pf}\left[ 
\begin{pmatrix}
0         & M_{01}    & M_{02}    & \dots & 1        & M_{0 2j+1} \\
-M_{01}   & 0         & M_{12}    & \dots & x        & M_{1 2j+1}  \\
\vdots    & \vdots    & \vdots    &       & \vdots   & \vdots \\
-1        & -x        & -x^2      & \dots & 0        & x^{2j+1} \\ 
-M_{0 2j+1} & -M_{1 2j+1} & -M_{2 2j+1} & \dots & -x^{2j+1}& 0 
\end{pmatrix}\right] .
\end{equation*}

To compute the skew-orthogonal polynomials numerically it is helpful to use
the following representation.  
There is a generalized Gaussian
elimination algorithm which factors the skew-symmetric matrix $M$ as
\begin{equation*}
M = L D^{1/2} J D^{1/2} L^T,
\end{equation*}
where $L$ is lower $2\times 2$ block triangular matrix
 with the $2\times 2$ block identity on the
diagonal, $J$ is given in (\ref{J}), 
and $D$ is a diagonal matrix with multiples of the $2\times 2$ block
identity on the diagonal.  
The rows of $L^{-1}$ give the coefficients of the skew-orthogonal
polynomials.
Numerical computation of skew-orthogonal polynomials with $V(x)
 =\frac{1}{2} x^2 + t x^4$ lead to
the conjecture 
 that the zeros of the even polynomials are real, interlace,
and their distribution converges to some measure; 
and that the zeros of the odd polynomials are real,
interlace, and their distribution converges to some measure.

The leading order of the asymptotics of the skew-orthogonal
polynomials have been computed by \cite{Eynard}.  
In section \ref{rh_sec} we
will compute a local $(d+1)\times (d+1)$ Riemann-Hilbert problem, where
$\mbox{deg}(V) = d$  (similar to that of
\cite{KM}) for the skew-orthogonal polynomials with respect to the
skew-inner product (\ref{GOE}).  
This Riemann-Hilbert problem may be tractable, at least in the sense
that the bi-orthogonal one is, and allow one to rigorously compute the
asymptotic expansion of the skew-orthogonal polynomials.  
The problem we find determines $p_{2k}$ uniquely but only finds $p_{2k+1}$
up to a multiple of $p_{2k}$.
An equivalent Riemann-Hilbert problem for the skew-orthogonal polynomials with
respect to the skew-inner product (\ref{GSE}) has not been found.

\section{\label{rh_sec} Riemann-Hilbert problem}

We will now work with the inner product (\ref{GOE}) however outside of
the random matrix context we may use the simpler expression 
\begin{equation}\label{skewinner}
\langle f, g \rangle_1 = 
\iint_{\mathbb{R}^2} f(x) g(y) \epsilon(x-y) e^{-V(x) - V(y)} dx dy,
\end{equation}
with no loss of generality. 
Assume that $V$ is a polynomial of degree $d$.  

Define 
\begin{equation*} %\label{wj}
w_j(x) = \int_\mathbb{R} y^j \epsilon(x-y) e^{-V(x)-V(y)} dy
\end{equation*}
and 
\begin{equation*} %\label{W}
\mathcal{W}(x) = e^{-2V(x)}.
\end{equation*}
Let 
\begin{equation*}
\langle f, g \rangle_2 = \int_{\mathbb{R}} f(x) g(x) \mathcal{W}(x) dx .
\end{equation*}

The skew inner product (\ref{skewinner}) 
is non-degenerate in the sense that 
the matrix $M_{ij} = \langle x^i, y^j \rangle_1$ is skew
diagonalizable; that is it can be written as 
\begin{equation*}
M = L J L^T
\end{equation*}
where $L$ is lower $2\times 2$ block triangular with non-zero
multiples of the $2\times 2$ identity matrix on the diagonal, and $J$
is given by (\ref{J}).  
An equivalent condition is that the principle $2n \times 2n$ minors of
$M$ are non-singular.  
The Pfaffian
 of these principle minors is proportional to a
multi-integral of a positive function, hence is non-zero.  
This non-degeneracy gives the
existence and uniqueness of the skew-orthogonal polynomials.

The family of skew-orthogonal polynomials with respect to
 (\ref{skewinner}) is
characterized by the conditions: 
$
\langle p_{2k}(x), y^j \rangle_1 = 0
$
and
$
\langle p_{2k+1}(x), y^j \rangle_1 = 0
$
for $0 \leq j \leq 2k-1$,  $p_{2k}(x) = x^{2k} +
\mathcal{O}(x^{2k-1}) $, and $p_{2k+1}(x) = x^{2k+1} +
\mathcal{O}(x^{2k-1}) $.

Define 
\begin{equation*}
\pi_{j+d-1}(y) = \frac{d}{dy} \left( y^j e^{-V(y)}\right) e^{V(y)} .
\end{equation*}
This function is a polynomial of degree $j+d-1$ in $y$.
The fundamental theorem of calculus and the definition of $\epsilon(x-y)$ 
implies that 
\begin{equation} \label{eq2.1}
\langle f(x), \pi_{j+d-1}(y) \rangle_1 = 2 \langle f(x), x^j \rangle_2.
\end{equation}
Therefore we find that the $2k$ orthogonality conditions on
$p_{2k}(x)$ become:  the $d-1$ conditions 
\begin{equation} \label{ortho1}
\langle p_{2k}(x), y^j \rangle_1 = 0, \; 0\leq j \leq d-2 
\end{equation}
together with the $2k-d+1$ conditions
\begin{equation} \label{ortho2}
\langle p_{2k}(x),  x^j \rangle_2 = 0, \; 0 \leq j \leq 2k-d.
\end{equation}
For $p_{2k+1}(x)$ we find the same conditions (with the same ranges on
$j$).

Consider a pair of Riemann-Hilbert problems:
Find a $(d+1) \times (d+1)$ matrix valued function $Y(z)$ satisfying:
\begin{enumerate}[1.]
\item Y is analytic on $\mathbb{C}/\mathbb{R}$.

\item Y has boundary values 
\begin{equation*}
Y_\pm(x) = \lim_{\epsilon \to 0^{\pm}} Y(x + i \epsilon)
\end{equation*}
for $x\in \mathbb{R}$

\item The boundary values of $Y$ satisfy the matrix equation 
\begin{equation} \label{rh}
Y_+ = Y_- M 
\end{equation}
where
\begin{equation*}
M = \begin{pmatrix} 1 & \mathcal{W}(x) & w_0(x) & \dots & w_{d-2}(x)
  \\
0 & 1 & 0 & \dots & 0  \\
\vdots &  & \ddots & \vdots \\
0 &   & \dots &  & 1  
\end{pmatrix}.
\end{equation*}

\item 
\begin{enumerate}[(a)]

\item
In the Even problem $Y$ satisfies the asymptotic boundary value
\begin{equation} \label{Even}
Y(z) = \left( I + \mathcal{O}(z^{-1})\right) \begin{pmatrix} 
z^{2k} &             &        & & \\
       & z^{-2k+d-1} &        & & \\
       &             & z^{-1} & & \\
       &             &        & \ddots & \\
       &             &        &        & z^{-1} 
\end{pmatrix} 
\end{equation}
as $ |z|\to \infty$.  

\item
In the Odd problem $Y$ satisfies the asymptotic boundary value
\begin{equation}\label{Odd}
Y(z) = \left( I + \mathcal{O}(z^{-1})\right)\begin{pmatrix} 
z^{2k+1} &             &        & & \\
       & z^{-2k+d-1} &        & & \\
       &             & z^{-1} & & \\
       &             &        & \ddots & \\
       &             &        &        & z^{-1} 
\end{pmatrix} 
\end{equation}
as $|z|\to \infty$.

\end{enumerate}
\end{enumerate}

 Let 
\begin{equation*} %\label{Cauchy}
C(f)(z) = \frac{1}{2\pi i} \int_\mathbb{R} \frac{f(x)}{x-z} dx
\end{equation*}
be the Cauchy Transform of $f(x)$.

We find two theorems:
\begin{thm}[Even] \label{eventhm}
The Riemann-Hilbert Problem above with asymptotic condition
(\ref{Even}) has a unique solution given by 
 the $(d+1) \times (d+1)$ matrix
\begin{equation*}
Y_{2k}(z) = \begin{pmatrix} 
p_{2k} & C( p_{2k} \mathcal{W}) & C( p_{2k} w_0 ) & \dots & C(p_{2k}
w_{d-2}) \\
p_{2k-1}^{(0)} &  \dots & &  \dots & \\
p_{2k-1}^{(1)} & \ddots & & \dots & \\
\vdots & \ddots & & & \\
p_{2k-1}^{(d-1)} & \dots & & \dots& 
\end{pmatrix}
\end{equation*}
where $p_{2k-1}^{(0)} = \alpha_k p_{2k-2} $ and 
$ p_{2k-1}^{(m)}, \; m>0$ are polynomials of degree at most $2k-1$
which will be specified in the proof.
The other entries of $Y_{2k}(z)$ are Cauchy transforms of
$p_{2k-1}^{(m)} \mathcal{W}$, and $p_{2k-1}^{(m)} w_n$, $0\leq n < d-1$. 
\end{thm} 

\begin{thm}[Odd] \label{oddthm}
The Riemann-Hilbert Problem above with asymptotic condition
(\ref{Odd})  has general solution given by 
the $(d+1) \times (d+1)$ matrix
\begin{equation*}
Y_{2k+1}(z) = \begin{pmatrix}
\tilde{p}_{2k+1} & C( \tilde{p}_{2k+1} \mathcal{W}) & 
C( \tilde{p}_{2k+1} w_0) & \dots & C(
\tilde{p}_{2k+1} w_{d-2}) \\ 
p_{2k}^{(0)} & \dots  & &  \dots & \\
p_{2k}^{(1)} & \ddots & & \dots & \\
\vdots & \ddots & & & \\
p_{2k}^{(d-1)} & \dots & & \dots & 
\end{pmatrix}
\end{equation*}
where $\tilde{p}_{2k+1} = p_{2k+1} + a_{k} p_{2k} $, 
$p_{2k}^{(0)} = b_{k, 0} p_{2k} +  c_{k, 0} p_{2k-2} $, 
and 
$p_{2k}^{(m)} = b_{k, m} p_{2k} + c_{k, m} p_{2k-1}^{(m)} $. 
The other entries of $Y_{2k+1}(z)$ are Cauchy transforms of $p_{2k}^{(m)}
\mathcal{W}$, and $p_{2k}^{(m)} w_n$, $0\leq n <d-1$.  In particular
 the general solution is determined up to a multiple of
 $p_{2k}$ in each row.  
\end{thm}

\section{Proof of theorems \ref{eventhm} and \ref{oddthm}}

Uniqueness of the solution in theorem \ref{eventhm} follows in the standard
way, by checking that both the jump condition (\ref{rh}) and the asymptotic
condition (\ref{Even})
have determinant one.  

To prove that $Y_{2k}$ is a solution of the Riemann-Hilbert problem with
asymptotic condition (\ref{Even}), one
checks that:
\begin{itemize}

\item $Y_{11}$ must be an analytic function of degree $2k$, hence a
  polynomial of degree $2k$.

\item $Y_{12}$ satisfies the jump condition
\begin{equation*}
Y_{12, +} = Y_{11} \mathcal{W}(x) + Y_{12, -}
\end{equation*}
across $\mathbb{R}$.
The Plemelj formula \cite{Ablowitz_Fokas} implies that 
\begin{equation*}
Y_{12}(z) = C( Y_{11} \mathcal{W} )(z).
\end{equation*}
For this function to have the asymptotics
\begin{equation*}
Y_{12}(z) = \mathcal{O}(z^{-2k+d-2}) 
\end{equation*}
the $2k-d+1$ conditions in (\ref{ortho2}) must be satisfied.

\item $Y_{1(j+3)}$ (for $0<j<d-2$) satisfies the jump condition
\begin{equation*}
Y_{1(j+3), +} = Y_{11} w_j(x) + Y_{1(j+3), -}
\end{equation*}
across $\mathbb{R}$.
The Plemelj formula implies that
\begin{equation*}
Y_{1(j+3)}(z) = C( Y_{11} w_j )(z).
\end{equation*}
For this function to have the asymptotics
\begin{equation*}
Y_{1(j+3)}(z) = \mathcal{O}(z^{-2})
\end{equation*} 
the condition (\ref{ortho1}) must be satisfied.

\end{itemize}

Hence $Y_{11}(z)$ is the polynomial $p_{2k}(z)$ satisfying conditions
(\ref{ortho1}) and (\ref{ortho2}).  

\subsection{Existence of the $p_{2k-1}^{(m)}$}  

The same arguments as above show that the lower rows of $Y_{2k}(z)$
are of the form 
\begin{equation} \label{Yrow}
\left( p_{2k-1}^{(m)}, C( p_{2k-1}^{(m)} \mathcal{W}), C( p_{2k-1}^{(m)}
  w_0 ) ,  \dots, C(p_{2k-1}^{(m)} w_{d-2} ) \right).
\end{equation}

Consider the second row:
The asymptotic conditions on $Y_{2k}(z)$ imply that $p_{2k-1}^{(0)}$
satisfies the following conditions as $|z|\to \infty$:
\begin{equation} \label{polycond}
p_{2k-1}^{(0)}(z) = \mathcal{O}(z^{2k-1}), 
\end{equation}
\begin{equation} \label{orthocond} 
C(p_{2k-1}^{(0)} \mathcal{W}) = 
- \frac{1}{2\pi i} \sum_{j=0}^\infty \frac{1}{z^{j+1}} \langle
p_{2k-1}^{(0)}, x^j \rangle_2 = z^{-2k+d-1} + \mathcal{O}(z^{-2k+d-2}) ,
\end{equation}
\begin{equation} \label{wmcond}
C(p_{2k-1}^{(0)} w_n ) = - \frac{1}{2\pi i} \sum_{j=0}^\infty
\frac{1}{z^{j+1}} \int_{\mathbb{R}} p_{2k-1}^{(0)} w_n(x) x^j dx =
  \mathcal{O}(z^{-2}), \; 0\leq n \leq d-2. 
\end{equation}
Condition (\ref{polycond}) forces $p_{2k-1}^{(0)}$ to be a polynomial
with degree at most $2k-1$.  Condition (\ref{orthocond}) is equivalent
to 
\begin{equation} \label{neworthocond}
 \langle p_{2k-1}^{(0)}, x^j \rangle_2 = \frac{1}{2} \langle 
p_{2k-1}^{(0)}(x),
  \pi_{j+d-1}(y) \rangle_1 = 0, \; 0\leq j \leq 2k-d-1 .
\end{equation}
The first $2k-d-1$ equations in (\ref{neworthocond}) together with the
$d-1$ conditions (\ref{wmcond}) imply that 
\begin{equation*}
p_{2k-1}^{(0)} = a_{2, 0}^{(2k-1)} p_{2k-1} + a_{1, 0}^{(2k-1)} p_{2k-2}.
\end{equation*}
The case $j=2k-d-1$ in (\ref{neworthocond}) implies that 
\begin{align*}
0 &= \langle p_{2k-1}^{(0)}(x), \pi_{2k-2}(y) \rangle_1 \\
&= a_{2, 0}^{(2k-1)} \langle
p_{2k-1}(x), \pi_{2k-2}(y) \rangle_1 + a_{1,0}^{(2k-1)} \langle p_{2k-2}(x),
\pi_{2k-2}(y) \rangle_1 
\\
&= a_{2,0}^{(2k-1)} \langle p_{2k-1}(x), y^{2k-2}\rangle_1, 
\end{align*}
from which we conclude that $a_{2,0}^{(2k-1)} = 0$.  
Finally one uses the leading order of (\ref{orthocond}) to find that
\begin{equation*}
a_{1,0}^{(2k-1)} = - 4\pi i \langle p_{2k-2}, y^{2k-1}\rangle_1^{-1} .
\end{equation*}
The non-degeneracy of the skew-inner product (\ref{skewinner}) guarantees
that $\langle p_{2k-2}, y^{2k-1}\rangle_1 \neq 0$.

Consider the other rows:
The asymptotic conditions on $Y_{2k}(z)$ imply that $p_{2k-1}^{(m)}, \; m>0$
satisfies the following conditions as $|z|\to \infty$:
\begin{equation} \label{mpolycond}
p_{2k-1}^{(m)}(z) = \mathcal{O}(z^{2k-1}), 
\end{equation}
\begin{equation} \label{morthocond} 
C(p_{2k-1}^{(m)} \mathcal{W}) = 
- \frac{1}{2\pi i} \sum_{j=0}^\infty \frac{1}{z^{j+1}} \langle
p_{2k-1}^{(m)}, x^j \rangle_2 = \mathcal{O}(z^{-2k+d-2}) ,
\end{equation}
\begin{equation} \label{mwmcond}
C(p_{2k-1}^{(m)} w_n ) = - \frac{1}{2\pi i} \sum_{j=0}^\infty
\frac{1}{z^{j+1}} \int_{\mathbb{R}} p_{2k-1}^{(m)} w_n(x) x^j dx =
\left\{ \begin{matrix}
  \mathcal{O}(z^{-2}) & m\neq n \\
z^{-1} + \mathcal{O}(z^{-2}) & m = n \end{matrix} \right. 
, \; 0\leq n \leq d-2. 
\end{equation}
Condition (\ref{mpolycond}) forces $p_{2k-1}^{(m)}$ to be a polynomial
with degree at most $2k-1$.
Condition (\ref{morthocond}) is equivalent
to 
\begin{equation} \label{mneworthocond}
 \langle p_{2k-1}^{(m)}, x^j \rangle_2=0, \; 0\leq j \leq 2k-d .
\end{equation}
The $2k-d+1$ conditions from (\ref{mneworthocond}) and the $d-2$
conditions from (\ref{mwmcond}) when $m\neq n$ determine the $2k-1$
free variables in $p_{2k-1}^{(m)}$ by solving a homogeneous linear
problem.  
One only need check that 
\begin{equation*}
\int_{\mathbb{R}} p_{2k-1}^{(m)} w_m(x) dx \neq 0.
\end{equation*}
Suppose that this integral is zero, then $\hat{p}_{2k} = p_{2k} +
p_{2k-1}^{(m)} $ satisfies the same orthogonality conditions as the
skew-orthogonal polynomial $p_{2k}$, which contradicts the uniqueness
of $p_{2k}$.  
This concludes the proof of Theorem \ref{eventhm}.

In fact we can say more: using (\ref{morthocond}) we see that we could
write the $p_{2k-1}^{(m)}$ in terms of $d-1$ of the orthogonal
polynomials with respect to $\mathcal{W}(x)$ given by (\ref{orthopoly}).

%%%%%%%%%%%%%%%%%%%%%%%%%%%%%%%%%%%%%%%%%%%%%%%%%%%%%%%%%%%%%%%%%%%%%%%
\subsection{The proof of Theorem \ref{oddthm}}

For the odd Riemann-Hilbert problem we note that
the asymptotic condition (\ref{Odd}) has determinant $z$.  So solutions will
only be unique up to possibly $d+1$ many degrees of freedom.  

The proof that the first row of $Y_{2k+1}$ has the form
\begin{equation*}
\left( \tilde{p}_{2k+1}(x), C(\tilde{p}_{2k+1} \mathcal{W}),
C(\tilde{p}_{2k+1} w_0), \dots, C(\tilde{p}_{2k+1} w_{d-2}) \right),
\end{equation*} 
follows as above.  The asymptotics of the first row of $Y_{2k}$ fit
within those of $Y_{2k+1}$ so we only get the first row of $Y_{2k+1}$
up to a multiple of the first row of $Y_{2k}$.  Therefore we may write
$ \tilde{p}_{2k+1} = p_{2k+1} + a_k p_{2k}$.  

The lower rows of $Y_{2k+1}$ have the form
\begin{equation*}
\left( p_{2k}^{(m)}, C(p_{2k}^{(m)} \mathcal{W}), C( p_{2k}^{(m)}
w_0), \dots, C(p_{2k}^{(m)} w_{d-2} ) \right).
\end{equation*}
We will find that the first entry of each row is only determined up to
a multiple of $p_{2k}$, as the orthogonality conditions on $p_{2k}$
fit within those of $p_{2k}^{(m)}$.  This gives $(d+1)$ free
parameters in the general solution to the Riemann-Hilbert problem.  

Consider the second row:
The asymptotic conditions on $Y_{2k+1}(z)$ imply that $p_{2k}^{(0)}$
satisfies the following conditions as $|z|\to \infty$:
\begin{equation} \label{oddpoly}
p_{2k}^{(0)} = \mathcal{O}(z^{2k}),
\end{equation}
\begin{equation} \label{oddorthocond}
C(p_{2k}^{(0)} \mathcal{W}) = - \frac{1}{2\pi i} \sum_{j=0}^\infty
\frac{1}{z^{j+1}} \langle p_{2k}^{(0)}, x^j\rangle_2 = z^{-2k+d-1} +
\mathcal{O}(z^{-2k+d-2}), 
\end{equation}
\begin{equation} \label{oddwncond}
C(p_{2k}^{(0)} w_n) = -\frac{1}{2\pi i} \sum_{j=0}^\infty
\frac{1}{z^{j+1}} \int_{\mathbb{R}} p_{2k}^{(0)} w_n(x) x^j dx =
\mathcal{O}(z^{-2}), \; 0\leq n \leq d-2.
\end{equation}
Condition (\ref{oddpoly}) forces $p_{2k}^{(0)}$ to be a polynomial
with degree at most $2k$.  Condition (\ref{oddorthocond}) is
equivalent to 
\begin{equation} \label{oddneworthocond}
 \langle p_{2k}^{(0)}, x^j \rangle_2 = \frac{1}{2} \langle p_{2k}^{(0)}(x),
\pi_{j+d-1}(y) \rangle_1=0, \; 0\leq j \leq 2k-d-1. 
\end{equation}
The first $2k-d-1$ equations in (\ref{oddneworthocond}) together with
the $d-1$ conditions (\ref{oddwncond}) imply that 
\begin{equation*}
p_{2k}^{(0)} = b_{k, 0} p_{2k} + d_{k,0} p_{2k-1} +
c_{k,0} p_{2k-2} .
\end{equation*}
The case $j=2k-d-1$ in (\ref{oddneworthocond}) implies that 
\begin{align*}
0 &= \langle p_{2k}^{(0)}(x), \pi_{2k-2}(y) \rangle_1 \\
&= b_{k,0} \langle p_{2k},
\pi_{2k-2}(y)\rangle_1 + d_{k,0} \langle p_{2k-1},
\pi_{2k-2}(y) \rangle_1  
+ c_{k,0} \langle p_{2k-2}, \pi_{2k-2}(y) \rangle_1 
\\
&= d_{k,0} \langle p_{2k-1}, y^{2k-2}\rangle_1 ,
\end{align*}
from which we conclude that $d_{k,0} = 0$.  Finally one uses the
leading order of (\ref{oddorthocond}) to  choose
 $c_{k,0}$ so that 
\begin{equation*}
c_{k,0}
\langle  p_{2k-2},
y^{2k-1}\rangle_1 = - 4\pi i.
\end{equation*}
Non-degeneracy of the inner product (\ref{skewinner}) guarantees that
this can be done.  The $b_{k, 0}$ remains as a free parameter.

Consider the other rows of $Y_{2k+1}$:  The asymptotic conditions on
$Y_{2k+1}(z)$ imply that $p_{2k}^{(m)}, \; m>0$ satisfy the following
conditions as $|z|\to \infty$:
\begin{equation} \label{oddmpolycond}
p_{2k}^{(m)} = \mathcal{O}(z^{2k}), 
\end{equation}
\begin{equation} \label{oddmorthocond}
C(p_{2k}^{(m)} \mathcal{W} ) = - \frac{1}{2\pi i} \sum_{j=0}^{\infty}
\frac{1}{z^{j+1}} \langle p_{2k}^{(m)}, x^j\rangle_2 
= \mathcal{O}(z^{-2k+d-2}),
\end{equation}
\begin{equation} \label{oddmwncond}
C(p_{2k}^{(m)} w_n) = - \frac{1}{2\pi i} \sum_{j=0}^\infty
\frac{1}{z^{j+1}} \int_{\mathbb{R}} p_{2k}^{(m)} w_n(x) x^j dx = 
\left\{ \begin{matrix} 
\mathcal{O}(z^{-2}) & m\neq n \\
z^{-1} + \mathcal{O}(z^{-2}) & m = n 
\end{matrix} \right.
, \; 0\leq n \leq d-2 
\end{equation}
Condition (\ref{oddmpolycond}) forces $p_{2k}^{(m)}$ to be a
polynomial with degree at most $2k$.  
The conditions (\ref{oddmpolycond}-\ref{oddmwncond}) are satisfied
 by $p_{2k-1}^{(m)}$, which fixes $2k-1$ of the free variables in
 $p_{2k}^{(m)}$.  
The orthogonality conditions satisfied by
$p_{2k}$ fit within those of $p_{2k}^{(m)}$, hence we find that
 $p_{2k}^{(m)} = b_{k, m} p_{2k} + c_{k, m} p_{2k-1}^{(m)}$ with
  $c_{k, m}$ determined by (\ref{oddmwncond}) with
 $n=m$, and $b_{k,m}$ a free parameter.  
This concludes the proof of Theorem \ref{oddthm}.

\section{Conclusion}

Both the bi-orthogonal and skew-orthogonal polynomials naturally split into a
pair of families: the bi-orthogonal polynomials from their definition, the
skew-orthogonal polynomials into the families of even and odd degree
polynomials.  
Numerical computations show that in both cases the distribution of zeros of
the two families appear to converge to
separate (yet mutually dependent) measures.  

The bi-orthogonal and skew-orthogonal polynomials are interesting
generalizations of the classical theory of orthogonal polynomials.  They
appear in the generalizations of random matrix theory from the well studied
GUE case.  So far little is understood about the necessary asymptotics for
both types.  We have formulated a Riemann-Hilbert problem for the
skew-orthogonal polynomials in the GOE case.  The structure of this
Riemann-Hilbert problem is close to that of \cite{KM} for the bi-orthogonal
polynomials.  An obvious problem with our formulation is that $p_{2k+1}$ is
not uniquely determined.  Experience with skew-orthogonal polynomials and the
Pfaff lattice hierarchy shows that this is a standard ambiguity complicating
the problems. 

The corresponding Riemann-Hilbert problem for the skew-orthogonal
polynomials associated to the GSE model has not been found.  The
equivalent to formula (\ref{eq2.1}) does not line up with the Cauchy
transform as nicely as in the $\beta=1$ case.

Our conclusions are: High order Riemann-Hilbert problems warrant further
investigation with the goal of finding methods to compute their asymptotics.  
In all studies of bi-orthogonal polynomials one should consider
simultaneously the equivalent result for skew-orthogonal polynomials.  One
suspects that there is a connection between these two types of polynomials.

\medskip

\noindent {\bf Acknowledgments:} I would like to thank
R. Buckingham and
K.T.-R. McLaughlin for many useful discussions and 
their comments and corrections.

\end{document}